\documentclass[reqno]{amsart}
\usepackage{bbold}
\usepackage{amsmath}
\usepackage{amssymb}
\usepackage{amsfonts}

\newcommand{\Abs}[1]{\left\vert #1 \right\vert}
\newcommand{\abs}[1]{\vert #1 \vert}

\newcommand{\norm}[1]{\left\Vert #1 \right\Vert}
\newcommand{\fixednorm}[1]{\Vert #1 \Vert}

\newcommand{\C}{\mathbb{C}}
\newcommand{\N}{\mathbb{N}}

\newcommand{\R}{\mathbb{R}}

\newcommand{\innerprod}[2]{\langle \, #1 , #2 \, \rangle}

\newcommand{\angles}[1]{\left\langle \, #1 \, \right\rangle}

\newtheorem{theorem}{Theorem}

\theoremstyle{definition}

\theoremstyle{remark}

\title{Ill-posedness of the Thirring model below the critical regularity}

\author[S.~Selberg]{Sigmund Selberg}
\author[A.~Tesfahun]{Achenef Tesfahun}

\address{Department of Mathematics, University of Bergen, PO Box 7803, 5020 Bergen, Norway}
\email{Sigmund.Selberg@uib.no}

\email{Achenef.Temesgen@uib.no}


\thanks{This work was partially supported by the project Pure Mathematics in Norway, funded by the Trond Mohn Foundation}

\begin{document}

\begin{abstract}
We consider a nonlinear $L^2$-critical nonlinear Dirac equation in one space dimension known as the Thirring model. Global well-posedness in $L^2$ for this equation was proved by Candy. Here we prove that the equation is ill posed in $L^p$ for $1 \le p < 2$, and in the massless case also in $H^s$ with $s < 0$.
\end{abstract}

\maketitle

\section{Introduction}

We consider the following nonlinear Dirac equation on the Minkowski space-time $\R^{1+1}$:
\begin{equation}\label{Thirring}
  (-i\gamma^\mu \partial_\mu + m) \psi = (\overline \psi \gamma^\mu \psi) \gamma_\mu \psi.
\end{equation}
In quantum field theory, this is known as the Thirring model \cite{Thirring} and describes the self-interactions of a Dirac field. Here we are interested in the Cauchy problem with initial data
\begin{equation}\label{SpinorData}
  \psi\vert_{t=0} = \psi_0
\end{equation}
for the Dirac spinor field $\psi \colon \R^{1+1} \to \C^2$, which is regarded as a column vector.

The equation is written in covariant form on $\R^{1+1}$ with coordinates $x^\mu$ ($\mu=0,1$) and metric $(g^{\mu\nu})=\mathrm{diag}(1,-1)$, where $x^0=t$ is time and $x^1=x$ is spatial position, and we write $\partial_\mu = \partial/\partial x^\mu$, so $\partial_0=\partial_t$ and $\partial_1=\partial_x$. The constant $m \ge 0$ is the particle mass, $\overline \psi = \psi^*\gamma^0$ is the adjoint spinor, where $\psi^*$ denotes the complex conjugate, and the $2 \times 2$ Dirac matrices $\gamma^0$ and $\gamma^1$ are required to satisfy
\[
  \gamma^\mu \gamma^\nu + \gamma^\nu \gamma^\mu = 2g^{\mu\nu} I
\]
as well as $(\gamma^0)^* = \gamma^0$ and $(\gamma^1)^* = -\gamma^1$.

Two key features of the Thirring model are the conservation of charge,
\begin{equation}\label{GlobalChargeConservation}
  \int_{\R} \abs{\psi(x,t)}^2 \, dx = \int_{\R} \abs{\psi_0(x)}^2 \, dx,
\end{equation}
and, in the case $m=0$, the invariance of \eqref{Thirring} with respect to the rescaling
\begin{equation}\label{ScalingLaw}
  \psi(x,t) \longrightarrow \psi^{(\lambda)}(x,t) := \lambda^{1/2} \psi(\lambda x,\lambda t) \qquad (\lambda > 0),
\end{equation}
which for initial data in $L^p(\R)$ gives
\[
  \fixednorm{\psi^{(\lambda)}(\cdot,0)}_{L^p(\R)} = \lambda^{1/2-1/p} \norm{\psi_0}_{L^p(\R)}.
\]
Thus, the data space $L^2(\R)$ (the charge space) is scaling-critical. Global well-posedness in that space (for any $m \ge 0$) was proved by Candy \cite{Candy2011}, and in the case $m=0$ also by Huh \cite{Huh2011}. We restate the result here for convenience.

\begin{theorem}[\cite{Candy2011}]\label{Thm0}
For any $\psi_0 \in L^2(\R)$ there is a global solution $\psi \in C(\R;L^2)$ of the Thirring model \eqref{Thirring} with the data \eqref{SpinorData}. Moreover, for any $T > 0$, the data-to-solution map $\psi_0 \mapsto \psi$ is continuous from $L^2(\R)$ into $C([-T,T];L^2(\R))$, and the solution is unique in a subset of that space. Higher regularity persists, in the sense that if $\psi_0 \in H^s(\R)$ for some $s > 0$, then $\psi \in C(\R;H^s)$.
\end{theorem}

Here $H^s(\R) = (1-\partial_x^2)^{-s/2} L^2(\R)$ is the standard $L^2$ based Sobolev space. We remark that the above result implies the same result in $L^2_\mathrm{loc}(\R)$, by the finite speed of propagation.

The analogous result for the related Gross-Neveu model was obtained by Huh and Moon \cite{Huh2015}, and for the Thirring and Gross-Neveu models coupled to an electromagnetic field by one of the authors in \cite{Selberg2018}. The first global existence result for the Thirring model is due to Delgado \cite{Delgado1978}, for more regular data, namely in $H^1(\R) = W^{1,2}(\R)$. A scattering result for the massive Thirring model has been obtained by Candy and Lindblad \cite{Candy2018}, while orbital stability of Dirac solitons in $L^2(\R)$ was proved by Contreras et al.~\cite{Contreras2016}.

The purpose of this note is to prove that local well-posedness fails below the critical regularity $L^2(\R)$, namely in $L^p(\R)$ with $1 \le p < 2$, which has the supercritical scaling behaviour: the data norm of the rescaled field $\psi^{(\lambda)}$ tends to infinity as $\lambda \to 0$, while the lifespan scales like $T \to T/\lambda$, suggesting that local well-posedness should fail. In the massless case, we additionally show ill-posedness in $H^s(\R)$ with $s < 0$ (note that the corresponding homogeneous norm again has the supercritical scaling).

To be precise, by \emph{local well-posedness} in a Banach space $X_0$ of initial data (which for us will be $L^p(\R)$ or $H^s(\R)$), containing $X_0 \cap L^2(\R)$ as a dense subspace, we mean here the following: For any data $\psi_0 \in X_0$ there exist
\begin{enumerate}
\item a neighbourhood $\Omega$ of $\psi_0$ in $X_0$,
\item a time $T > 0$, and
\item a continuous map $S \colon \Omega \to C([0,T];X_0)$ which on $\Omega \cap L^2(\R)$ agrees with the data-to-solution map from Theorem \ref{Thm0}.
\end{enumerate}

This notion of local well-posedness is fairly weak, since we do not require that for every $\psi_0' \in \Omega$, the field $\psi' := S(\psi_0')$ must satisfy the nonlinear Dirac equation in the sense of distributions on $\R \times (0,T)$. All we require is that $S$ be an extension of the data-to-solution map on $\Omega \cap L^2(\R)$, which exists by Theorem \ref{Thm0}.

Our main result is that even this rather weak notion of local well-posedness fails for the Thirring model when $X_0=L^p(\R)$ with $p < 2$.

\begin{theorem}\label{Thm1}
The Cauchy problem \eqref{Thirring}, \eqref{SpinorData} for the Thirring model fails to be locally well posed (in the above sense) in the data space $L^p(\R)$ for $1 \le p < 2$. Moreover, in the massless case $m=0$, local well-posedness fails also in $H^s(\R)$ for $s < 0$.
\end{theorem}

More precisely, the failure of local well-posedness is demonstrated for the specific choice of data
\[
  \psi_0(x) = \chi_{(-1,1)}(x) \frac{1}{\abs{x}^{1/2}} \begin{pmatrix} 1 \\ 1 \end{pmatrix},
\]
where $\chi_{(-1,1)}$ is the characteristic function of the interval $(-1,1)$. These data belong to $L^p(\R)$ for $1 \le p < 2$, but fail to belong to $L^2(\R)$. Moreover, by dual Sobolev embedding (or by showing directly that the Fourier transform is $O(\angles{\xi}^{-1/2})$) it follows that they also belong to $H^s(\R)$ for $s < 0$.

We then show that, regardless of the choice of neighbourhood $\Omega$ of $\psi_0$ in $L^p$, and regardless of how small we take $T > 0$, there cannot exist an extension $S : \Omega \to C([0,T];L^p)$ of the $L^2$-solution map from Theorem \ref{Thm0}. To this end, we consider the approximating data
\[
  \psi_0^\varepsilon(x) = \chi_{(-1,1)}(x) \frac{1}{(\varepsilon+\abs{x})^{1/2}} \begin{pmatrix} 1 \\ 1 \end{pmatrix}
\]
for $\varepsilon > 0$, which belong to $L^2(\R)$ and therefore evolve to globally well posed solutions $\psi^\varepsilon \in C(\R;L^2)$ by Theorem \ref{Thm0}. Note that as $\varepsilon \to 0$, $\psi^\varepsilon_0 \to \psi_0$ in $L^p(\R)$ for $1 \le p < 2$, and hence also in $H^s(\R)$ for $s < 0$.

So to disprove local well-posedness in $L^p$, it suffices to show that $\psi^\varepsilon$ does not have a limit in $C([0,T];L^p)$, no matter how small we take $T > 0$, and similarly in the case of $H^s$.

Before turning to the proof of Theorem \ref{Thm1}, we make some remarks regarding the relation of our choice of data to the scaling law \eqref{ScalingLaw}, in the massless case $m=0$. First, we note that the data cut-off $\chi_{(-1,1)}$ was put in simply to avoid having to work with local $L^p$ spaces, and on account of the finite speed of propagation it will not play any significant role in our analysis. So ignoring for now the cut-off, the data are then scale-invariant:
\begin{equation}\label{ScaleInvariantData}
  \psi_0(x) = \frac{1}{\abs{x}^{1/2}} \begin{pmatrix} 1 \\ 1 \end{pmatrix}.
\end{equation}
If there exists a unique solution corresponding to these data, it must necessarily be scale invariant (that is, $\psi^{(\lambda)}=\psi$ for all $\lambda > 0$) and hence self-similar:
\begin{equation}\label{Selfsimilar}
  \psi(x,t) = \frac{1}{\sqrt{t}} \phi\left(\frac{x}{t}\right) \qquad (t > 0)
\end{equation}
with
\begin{equation}\label{SelfsimilarData}
  \psi_0(x) = \lim_{t \to 0^+} \frac{1}{\sqrt{t}} \phi\left(\frac{x}{t}\right).
\end{equation}
This leads us to the interesting question whether there exist any such solutions (other than the trivial one---the zero solution). The answer is No, as we prove in Theorem \ref{Thm2} below. In particular, this shows that either the existence or the uniqueness of a solution for the scale invariant data \eqref{ScaleInvariantData} must fail. In fact, we shall see that it is in some sense the uniqueness that fails. Indeed, it can be observed (see section \ref{FinalRemarks}) that the approximating solutions $\psi^{\varepsilon}$ (without the cut-off on the data) converge to a continuum of possible solutions, depending on the sequence along which $\varepsilon > 0$ tends to zero.

We now state our result on non-existence of self-similar solution. For this, we use the following terminology: Suppose we are given a data space $X_0$, some data $\psi_0 \in X_{0,\mathrm{loc}}$, and a solution $\psi \in C([0,T];X_{0,\mathrm{loc}})$ with these data, for some $T > 0$. We say that the solution is \emph{stable} if, for any sequence of data $\psi_0^{(n)} \in X_{0,\mathrm{loc}} \cap L^2_\mathrm{loc}(\R)$ such that $\fixednorm{\psi_0^{(n)} - \psi_0}_{X_{0,\mathrm{loc}}} \to 0$ as $n \to \infty$, we have $\norm{\psi^{(n)}(t) - \psi(t)}_{X_{0,\mathrm{loc}}} \to 0$ as $n \to \infty$, for $t \in [0,T]$. Here $\psi^{(n)} \in C(\R;L^2_\mathrm{loc})$ is the solution with data $\psi_0^{(n)}$, which exists by Theorem \ref{Thm0}. 

We have the following result.

\begin{theorem}\label{Thm2}
Let $X_0$ be $L^p(\R)$ or $H^s$, for any $1 \le p \le \infty$ or $s \in \R$. Then there does not exist any stable, non-zero solution $\psi \in C([0,\infty);X_{0,\mathrm{loc}})$ of \eqref{Thirring} satisfying the self-similarity assumptions \eqref{Selfsimilar}, \eqref{SelfsimilarData}.
\end{theorem}

We remark that for $p \ge 2$, respectively $s \ge 0$, the conclusion follows from Theorem \ref{Thm0}, which of course excludes the self-similar blow-up at $(x,t)=(0,0)$. In this connection we mention also that non-existence of self-similar solutions for a class of nonlinear Dirac equations in one space dimension (which does not include the Thirring model, however) has been proved recently by Huh and Pelinovsky \cite{Huh2019}.

The paper is organized as follows: In the next section, we rewrite the system in terms of the components of the spinor, and we recall one of its key features: the conservation of charge (in its local form). In section \ref{Massless} we consider first the massless case, where the solution $\psi^\varepsilon$ can be computed directly and the ill-posedness is easily deduced. However, it is not clear how this argument can be generalised to the massive case. To get around this problem, we present in section \ref{Massless2} a different argument for the massless case, motivating the proof of the massive case, given in section \ref{Massive}. The non-existence of self-similar solutions is proved in section \ref{SelfsimilarSection}. Finally, in section \ref{FinalRemarks} we make some further remarks on the precise behaviour of the solutions in the massless case, observing that $\psi^\varepsilon$ bifurcates to a continuum of possible solutions depending on the sequence along which $\varepsilon$ is allowed to tend to zero. 

\section{Preliminaries}\label{Prelims}

Adopting the particular representation
\[
  \gamma^0 =
  \begin{pmatrix}
    0 & 1  \\
    1 &0
  \end{pmatrix},
  \qquad
  \gamma^1=
  \begin{pmatrix}
    0 & -1  \\
    1 & 0
  \end{pmatrix}
\]
for the Dirac matrices, and writing $\psi = (u,v)^\intercal$, where $u$ and $v$ take values in $\C$, we see that the Cauchy problem takes the convenient form
\begin{equation}\label{CP}
\left\{
\begin{alignedat}{2}
 (\partial_t + \partial_x)u &= -imv + 2i \abs{v}^2 u,& \qquad u(x,0) &= f(x),
 \\
 (\partial_t - \partial_x)v &= -imu + 2i \abs{u}^2 v,& \qquad v(x,0) &= g(x).
\end{alignedat}
\right.
\end{equation}

To get ill-posedness we use the data
\begin{equation}\label{Data}
  f(x) = g(x) = \frac{\chi_{(-1,1)}(x)}{\abs{x}^{1/2}}
\end{equation}
and the approximations
\begin{equation}\label{ApproxData}
  f_\varepsilon(x) = g_\varepsilon(x) = \frac{\chi_{(-1,1)}(x)}{(\varepsilon + \abs{x})^{1/2}} \qquad (\varepsilon > 0).
\end{equation}
The latter belong to $L^2(\R)$, hence there is a corresponding global solution to \eqref{CP} in $C(\R;L^2(\R))$, which we denote $(u_\varepsilon,v_\varepsilon)$. Note that $f_\varepsilon$ converges to $f$ in $L^p(\R)$ for $1 \le p < 2$, hence also in $H^s(\R)$ for $s < 0$.

In fact, we will remove the characteristic function $\chi_{(-1,1)}$ when performing our calculations. By finite speed of propagation, this does not affect the solution in the region $\abs{x} + \abs{t} < 1$, and for the proof of ill-posedness it will be enough to work in that region.

A key fact that we will make use of is the conservation law $\partial_\mu j^\mu = 0$, where $j^\mu = \overline{\psi}\gamma^\mu\psi$ is the Dirac charge density. Integrating this over the time-slab $\R \times [0,t]$ gives the conservation of global charge \eqref{GlobalChargeConservation}, while integrating over a backward cone gives the local form of conservation of charge,
\begin{equation}\label{LocalChargeConservation}
  \int_0^t 2\abs{u(x+t-\sigma,\sigma)}^2 \, d\sigma + \int_0^t 2\abs{v(x-t+\sigma,\sigma)}^2 \, d\sigma
  = \int_{x-t}^{x+t} \left( \abs{f(y)}^2 + \abs{g(y)}^2 \right) \, dy.
\end{equation}
See, e.g., section 2 of \cite{Selberg2018}.

\section{The massless case}\label{Massless}

In this section we assume $m=0$, in which case the system
\[
\begin{alignedat}{2}
 (\partial_t + \partial_x)u &= 2i \abs{v}^2 u,& \qquad u(x,0) &= f(x),
 \\
 (\partial_t - \partial_x)v &= 2i \abs{u}^2 v,& \qquad v(x,0) &= g(x)
\end{alignedat}
\]
can be integrated by multiplying $u$ and $v$ by integrating factors $e^{-i\phi_+}$ and $e^{-i\phi_-}$, respectively, where $\phi_+$ and $\phi_-$ are defined by
\[
\begin{alignedat}{2}
 (\partial_t + \partial_x)\phi_+ &= 2\abs{v}^2,& \qquad \phi_+(x,0) &= 0,
 \\
 (\partial_t - \partial_x)\phi_- &= 2\abs{u}^2,& \qquad \phi_-(x,0) &= 0,
\end{alignedat}
\]
that is,
\begin{align*}
  \phi_+(t,x) &= \int_0^t 2\abs{v(x-t+\sigma,\sigma)}^2 \, d\sigma,
  \\
  \phi_-(t,x) &= \int_0^t 2\abs{u(x+t-\sigma,\sigma)}^2 \, d\sigma.
\end{align*}
This gives
\begin{align*}
  u(x,t) &= f(x-t) e^{i\phi_+(t,x)},
  \\
  v(x,t) &= g(x+t) e^{i\phi_-(t,x)}.
\end{align*}
In particular, $\abs{u(x,t)} = \abs{f(x-t)}$ and $\abs{v(x,t)} = \abs{g(x+t)}$, so we have
\begin{align*}
  \phi_+(x,t) &= \int_0^t 2\abs{g(x-t+2\sigma)}^2 \, d\sigma
  = \int_{x-t}^{x+t} \abs{g(s)}^2 \, ds,
  \\
  \phi_-(x,t) &= \int_0^t 2\abs{f(x+t-2\sigma)}^2 \, d\sigma
  = \int_{x-t}^{x+t} \abs{f(s)}^2 \, ds.
\end{align*}

Now we choose our data as in \eqref{ApproxData}, but without the characteristic function (see the remark in section \ref{Prelims}) and calculate the solution $(u_\varepsilon,v_\varepsilon)$. Thus,
\[
  f_\varepsilon(x) = g_\varepsilon(x) = \frac{1}{(\varepsilon + \abs{x})^{1/2}}.
\]
Since $f_\varepsilon=g_\varepsilon$, we have $\phi_+^\varepsilon=\phi_-^\varepsilon$ and we denote this function simply by $\phi_\varepsilon$. It is given by
\[
  \phi_\varepsilon(x,t) = \int_{x-t}^{x+t} \frac{dy}{\varepsilon+\abs{y}},
\]
and we calculate it for $t > 0$ by separating into the regions corresponding to $x \ge t$, $t \ge \abs{x}$ and $x \le -t$, respectively:
\begin{enumerate}
  \item For $x \ge t$,
  \[
    \phi_\varepsilon(x,t) = \int_{x-t}^{x+t} \frac{dy}{\varepsilon+y} = \log(\varepsilon+x+t) - \log(\varepsilon+x-t);
  \]
  \item For $t \ge \abs{x}$,
  \[
    \phi_\varepsilon(x,t) = \int_{x-t}^{0} \frac{dy}{\varepsilon-y} + \int_0^{x+t} \frac{dy}{\varepsilon+y} = -2\log\varepsilon + \log(\varepsilon+t-x) + \log(\varepsilon+x+t);
  \]
  \item For $x \le -t$,
  \[
    \phi_\varepsilon(x,t) = \int_{x-t}^{x+t} \frac{dy}{\varepsilon-y} = \log(\varepsilon+t-x) - \log(\varepsilon-x-t).
  \]
\end{enumerate}

Thus, in the region $t > \abs{x}$,
\begin{align*}
  e^{2i\log\varepsilon} u_\varepsilon(x,t) &= \frac{1}{(\varepsilon+t-x)^{1/2}} e^{i[\log(\varepsilon + t-x) + \log(\varepsilon + t+x)]},
  \\
  e^{2i\log\varepsilon} v_\varepsilon(x,t) &= \frac{1}{(\varepsilon+t+x)^{1/2}} e^{i[\log(\varepsilon + t-x) + \log(\varepsilon + t+x)]},
\end{align*}
where the right hand sides converge uniformly (as do all derivatives) on compact subsets of $t > \abs{x}$, as $\varepsilon \to 0$. Fixing $t \in (0,1/2)$ and restricting to $x \in (-t,t)$, it follows that $(u_\varepsilon,v_\varepsilon)$ cannot converge in $L^p$ or in $H^s$, as $\varepsilon \to 0$. On the other hand, the data converge in those spaces if $p < 2$ or $s < 0$, respectively, so we have a contradiction to local well-posedness in those cases.

This concludes the proof of Theorem \ref{Thm1} in the massless case. The problem is that it is not clear how to generalize the above argument to the massive case. In the next section we provide an alternative argument for the ill-posedness in $L^p$ in the massless case, which is less direct but gives us the idea of how to treat the massive case.

\section{The massless case---an alternative approach}\label{Massless2}

The argument for the massless case in the previous section depends strongly on the fact that $\phi_+^\varepsilon$ and $\phi_-^\varepsilon$ are explicitly known, so that we can integrate the equations. In the massive case, on the other hand, these functions are not explicitly known, so it is hard to see how to generalize the argument. What one does know also in the massive case, however, is that
\[
  \phi_+^\varepsilon(x,t) + \phi_-^\varepsilon(x,t) = \int_{x-t}^{x+t} \frac{dy}{\varepsilon+\abs{y}},
\]
by the conservation of charge \eqref{LocalChargeConservation}. To make use of this fact, it is natural to try to deduce the ill-posedness in terms of the product $uv$, since then the above phase appears. We first do this in the massless case to show the idea in a simple setting, and then in the next section we show how to generalize to the massive case.

In the massless case,
\[
  e^{-i(\phi_+^\varepsilon + \phi_-^\varepsilon)} u_\varepsilon v_\varepsilon (x,t) = f_\varepsilon(x-t)g_\varepsilon(x+t).
\]
In the region $t > \abs{x}$, this becomes
\[
  e^{4i\log\varepsilon} u_\varepsilon v_\varepsilon (x,t)
  = \frac{e^{2i\log(\varepsilon+x+t)}}{(\varepsilon+x+t)^{1/2}}
  \frac{e^{2i\log(\varepsilon+t-x)}}{(\varepsilon+t-x)^{1/2}}.
\]
Choosing positive sequences $\varepsilon_n \to 0$ and $\varepsilon_n' \to 0$ such that $e^{4i\log\varepsilon_n} = 1$ and $e^{4i\log\varepsilon_n'} = -1$ for all $n \in \N$, we get
\begin{equation}\label{KeyFact1}
  + u_\varepsilon v_\varepsilon (x,t)
  = \frac{e^{2i\log(\varepsilon+x+t)}}{(\varepsilon+x+t)^{1/2}}
  \frac{e^{2i\log(\varepsilon+t-x)}}{(\varepsilon+t-x)^{1/2}} \quad \text{with $\varepsilon=\varepsilon_n$},
\end{equation}
as well as
\begin{equation}\label{KeyFact2}
  - u_\varepsilon v_\varepsilon (x,t)
  = \frac{e^{2i\log(\varepsilon+x+t)}}{(\varepsilon+x+t)^{1/2}}
  \frac{e^{2i\log(\varepsilon+t-x)}}{(\varepsilon+t-x)^{1/2}}
  \quad \text{with $\varepsilon=\varepsilon_n'$}.
\end{equation}
Note that as $n \to \infty$, the right hand sides converge pointwise in the region $t > \abs{x}$, to
\[
  \frac{e^{2i\log(x+t)}}{(x+t)^{1/2}} \frac{e^{2i\log(t-x)}}{(t-x)^{1/2}}.
\]

Now let $p \in [1,2)$ and assume that local well-posedness holds in $L^p$. Then, since our initial data $(f_\varepsilon,g_\varepsilon)$ converge in $L^p$ to $(f,g)$, the solution $(u_\varepsilon,v_\varepsilon)$ has a limit $(u,v)$ in $C([0,T];L^p)$ for $T > 0$ sufficiently small. But convergence in $L^p$ implies pointwise converges a.e.~of a subsequence. Thus, fixing a $t \in (0,T)$ and passing to suitable subsequences of $(\varepsilon_n)$ and $(\varepsilon_n')$, we may assume that $u_\varepsilon(x,t) \to u(x,t)$ and $v_\varepsilon(x,t) \to v(x,t)$, for almost every $x$, along $\varepsilon=\varepsilon_n$ and also along $\varepsilon = \varepsilon_n'$.

But then, passing to the limit in \eqref{KeyFact1} and \eqref{KeyFact2}, we get
\[
  + uv(x,t) = - uv(x,t) = \frac{e^{2i\log(x+t)}}{(x+t)^{1/2}} \frac{e^{2i\log(t-x)}}{(t-x)^{1/2}}
\]
for almost every $x \in (-t,t)$. But this is clearly not possible, since the right hand side is nonzero everywhere in that interval. The assumption of local well-posedness therefore leads to a contradiction.

We now show how to generalize this argument to the massive case.

\section{The massive case}\label{Massive}

In this section we consider the general case, $m \ge 0$, so the Cauchy problem reads
\[
\begin{alignedat}{2}
 (\partial_t + \partial_x)u &= -imv + 2i \abs{v}^2 u,& \qquad u(x,0) &= f(x),
 \\
 (\partial_t - \partial_x)v &= -imu + 2i \abs{u}^2 v,& \qquad v(x,0) &= g(x).
\end{alignedat}
\]
Again we define $\phi_+$ and $\phi_-$ by
\[
\begin{alignedat}{2}
 (\partial_t + \partial_x)\phi_+ &= 2\abs{v}^2,& \qquad \phi_+(x,0) &= 0,
 \\
 (\partial_t - \partial_x)\phi_- &= 2\abs{u}^2,& \qquad \phi_-(x,0) &= 0,
\end{alignedat}
\]
that is,
\begin{align*}
  \phi_+(x,t) &= \int_0^t 2\abs{v(x-t+\sigma,\sigma)}^2 \, d\sigma,
  \\
  \phi_-(x,t) &= \int_0^t 2\abs{u(x+t-\sigma,\sigma)}^2 \, d\sigma.
\end{align*}

Now
\[
\begin{alignedat}{2}
 (\partial_t + \partial_x)(e^{-i\phi_+}u) &= -ime^{-i\phi_+}v,& \qquad e^{-i\phi_+}u(0,x) &= f(x),
 \\
 (\partial_t - \partial_x)(e^{-i\phi_-}v) &= -ime^{-i\phi_-}u,& \qquad e^{-i\phi_-}v(0,x) &= g(x),
\end{alignedat}
\]
which integrates to
\begin{align}
  \label{IntegratedA}
  e^{-i\phi_+}u(x,t) &= f(x-t) - im \int_0^t e^{-i\phi_+}v (x-t+\sigma,\sigma) \, d\sigma,
  \\
  \label{IntegratedB}
  e^{-i\phi_-}v(x,t) &= g(x+t) - im \int_0^t e^{-i\phi_-}u (x+t-\sigma,\sigma) \, d\sigma.
\end{align}
Thus,
\[
  e^{-i(\phi_+ + \phi_-)}uv(x,t)
  =
  f(x-t)g(x+t) + \sum_{j=1}^3 R_j(x,t),
\]
where
\begin{align*}
  R_1(x,t) &= f(x-t) \left( - im \int_0^t e^{-i\phi_-}u (x+t-\sigma,\sigma) \, d\sigma \right),
  \\
  R_2(x,t) &= g(x+t) \left( - im \int_0^t e^{-i\phi_+}v (x-t+\sigma,\sigma) \, d\sigma \right),
  \\
  R_3(x,t) &= -m^2\left( \int_0^t e^{-i\phi_-}u (x+t-\sigma,\sigma) \, d\sigma \right) \left( \int_0^t e^{-i\phi_+}v (x-t+\sigma,\sigma) \, d\sigma \right).
\end{align*}

Now we take the data $(f_\varepsilon,g_\varepsilon)$ as in \eqref{ApproxData} and introduce superscripts or subscripts $\varepsilon$ on the functions to indicate their dependence on that parameter. By the conservation of charge \eqref{LocalChargeConservation},
\[
  (\phi_+^\varepsilon + \phi_-^\varepsilon)(x,t) = \int_{x-t}^{x+t} \frac{dy}{\varepsilon+\abs{y}},
\]
where the right hand side was calculated in section \ref{Massless}. Then for $t > \abs{x}$,
\begin{equation}\label{KeyFormula}
  e^{4i\log\varepsilon}u_\varepsilon v_\varepsilon(x,t)
  =
  \frac{e^{2i\log(\varepsilon+x+t)}}{(\varepsilon+x+t)^{1/2}}
  \frac{e^{2i\log(\varepsilon+t-x)}}{(\varepsilon+t-x)^{1/2}}
  + R_\varepsilon(x,t),
\end{equation}
where
\[
  R_\varepsilon(x,t) = e^{2i\log(\varepsilon+x+t)} e^{2i\log(\varepsilon+t-x)}\sum_{j=1}^3 R_{j,\varepsilon}(x,t).
\]

Thus, to make the argument from massless case in the previous section work out, all we have to do is to show that the remainder terms $R_{j,\varepsilon}$ are negligible compared to the first term on the right hand side of \eqref{KeyFormula}, pointwise in some closed ball in the region $t > \abs{x}$.

We choose the ball $B = B_{\delta/4}(0,\delta)$ centered at $(x,t) = (0,\delta)$ with radius $\delta/4$, where $\delta \ll 1$ remains to be chosen. Assuming $\varepsilon < \delta$ (as we may, since $\varepsilon$ will tend to zero), then on the absolute value of the first term on the right hand side of \eqref{KeyFormula} we have the bounds
\begin{equation}\label{MainBound}
  \frac{1}{2\delta}
  \le \Abs{ \frac{e^{2i\log(\varepsilon+x+t)}}{(\varepsilon+x+t)^{1/2}}
  \frac{e^{2i\log(\varepsilon+t-x)}}{(\varepsilon+t-x)^{1/2}} }
  \le \frac{2}{\delta}
  \qquad \text{for $(x,t) \in B$.}
\end{equation}

The next step is to derive some pointwise bounds on $\abs{R_{j,\varepsilon}}$. To this end, it suffices to obtain bounds for the integrals
\[
  \int_0^t \abs{u_\varepsilon(x+t-\sigma,\sigma)} \, d\sigma
  \quad \text{and} \quad
  \int_0^t \abs{v_\varepsilon(x-t+\sigma,\sigma)} \, d\sigma.
\]
To obtain such bounds we define, for $t > 0$,
\[
  A(t) = \sup_{y \in \R} \int_0^t \abs{u_\varepsilon(y-\sigma,\sigma)} \, d\sigma + \sup_{y \in \R} \int_0^t \abs{v_\varepsilon(y+\sigma,\sigma)} \, d\sigma.
\]
From \eqref{IntegratedA} and \eqref{IntegratedB} we have, for any $y \in \R$ and $\sigma > 0$,
\begin{align*}
  \abs{u_\varepsilon(y-\sigma,\sigma)} &\le \frac{1}{(\varepsilon+\abs{y-2\sigma})^{1/2}} + m \int_0^\sigma \abs{v_\varepsilon(y-2\sigma+s,s)} \, ds,
  \\
  \abs{v_\varepsilon(y+\sigma,\sigma)} &\le \frac{1}{(\varepsilon+\abs{y+2\sigma})^{1/2}} + m \int_0^\sigma \abs{u_\varepsilon(y+2\sigma-s,s)} \, ds,
\end{align*}
and integrating this with respect to $\sigma \in (0,t)$ we get
\begin{align*}
  \int_0^t \abs{u_\varepsilon(y-\sigma,\sigma)} \, d\sigma &\le \int_0^t \frac{d\sigma}{\abs{y-2\sigma}^{1/2}} + m \int_0^t \int_0^\sigma \abs{v_\varepsilon(y-2\sigma+s,s)} \, ds \, d\sigma,
  \\
  \int_0^t \abs{v_\varepsilon(y+\sigma,\sigma)} \, d\sigma &\le \int_0^t \frac{d\sigma}{\abs{y+2\sigma}^{1/2}} + m \int_0^t \int_0^\sigma \abs{u_\varepsilon(y+2\sigma-s,s)} \, ds \, d\sigma,
\end{align*}
hence
\begin{align*}
  \sup_{y \in \R} \int_0^t \abs{u_\varepsilon(y-\sigma,\sigma)} \, d\sigma &\le 4 t^{1/2} + m \int_0^t A(\sigma) \, d\sigma,
  \\
  \sup_{y \in \R} \int_0^t \abs{v_\varepsilon(y+\sigma,\sigma)} \, d\sigma &\le 4 t^{1/2} + m \int_0^t A(\sigma) \, d\sigma,
\end{align*}
implying
\[
  A(t) \le 8 t^{1/2} + 2m \int_0^t A(\sigma) \, d\sigma.
\]
By Gr\"onwall's inequality it follows that
\[
  A(t) \le 8 t^{1/2} e^{2mt}
\]
for $t > 0$.

We can now obtain the desired bounds for $\abs{R_{j,\varepsilon}}$ in the ball $B = B_{\delta/4}(0,\delta)$ for $\delta > 0$ sufficiently small. Clearly,
\[
  \sum_{j=1}^3 \abs{R_{j,\varepsilon}(x,t)} \le \frac{mA(t)}{\abs{x-t}^{1/2}} + \frac{mA(t)}{\abs{x+t}^{1/2}} + [mA(t)]^2,
\]
and for $(x,t) \in B$ we have $\abs{x \pm t} \ge \delta/2$ and $A(t) \le c \delta^{1/2}$ (where $c=16e^{2m}$ suffices), so we get an upper bound
\begin{equation}\label{RemainderBound}
  \sup_B \abs{R_{j,\varepsilon}} \le C,
\end{equation}
uniformly in $\delta$ and in $\varepsilon < \delta$ (in fact, we can take $C=4mc + m^2c^2$, where $c=16e^{2m}$).

Comparing \eqref{MainBound} and \eqref{RemainderBound}, we see that on the right hand side of \eqref{KeyFormula}, the first term dominates for $(x,t)$ in $B = B_{\delta/4}(0,\delta)$ if we choose $0 < \delta \ll 1/C$. It is then clear that we can carry through the argument from the massless case in the previous section and we get a contradiction to well-posedness.

This concludes the proof of Theorem \ref{Thm1}.

\section{Non-existence of self-similar solutions}\label{SelfsimilarSection}

Here we prove Theorem \ref{Thm2}. Let $X_0$ be $L^p(\R)$ or $H^s$, for any $1 \le p < 2$ or $s < 0$. Assume that there exists a stable, non-zero solution $\psi \in C([0,\infty);X_{0,\mathrm{loc}})$ of \eqref{Thirring} satisfying the self-similarity assumptions
\[
  \psi(x,t) = \frac{1}{\sqrt{t}} \phi\left(\frac{x}{t}\right) \qquad (t > 0)
\]
and
\[
  \psi_0(x) = \lim_{t \to 0^+} \frac{1}{\sqrt{t}} \phi\left(\frac{x}{t}\right) \quad \text{in $X_{0,\mathrm{loc}}$}.
\]
Then necessarily $\lambda^{1/2} \psi_0(\lambda x) = \psi_0(x)$, so the data are scale invariant and therefore of the form
\[
  \psi_0(x) = \begin{pmatrix} f(x) \\ g(x) \end{pmatrix},
  \qquad f(x) = \frac{\kappa_\pm}{\abs{x}^{1/2}},
  \qquad g(x) = \frac{\lambda_\pm}{\abs{x}^{1/2}}
  \qquad \text{for $\pm x > 0$},
\]
for some constants $\kappa_+,\kappa_-,\lambda_+,\lambda_- \in \C$. Let $f_\varepsilon, g_\varepsilon \in L^2_{\mathrm{loc}}(\R)$ for $\varepsilon > 0$ be the approximating data obtained by replacing the denominator $\abs{x}^{1/2}$ by $(\varepsilon+\abs{x})^{1/2}$. Then as in section \ref{Massless} we compute the corresponding solution $\psi^\varepsilon=(u_\varepsilon,v_\varepsilon)^\intercal$ to be, for $t > 0$:
\begin{enumerate}
  \item For $x \ge t$,
  \begin{align*}
    u_\varepsilon(x,t) &= \frac{\kappa_+}{(\varepsilon+x-t)^{1/2}} e^{i\abs{\lambda_+}^2[\log(\varepsilon+x+t) - \log(\varepsilon+x-t)]},
    \\
    v_\varepsilon(x,t) &= \frac{\lambda_+}{(\varepsilon+x+t)^{1/2}} e^{i\abs{\kappa_+}^2[\log(\varepsilon+x+t) - \log(\varepsilon+x-t)]};
  \end{align*}
  \item For $t \ge \abs{x}$,
  \begin{align*}
    u_\varepsilon(x,t) &= \frac{\kappa_-}{(\varepsilon+t-x)^{1/2}} 
    e^{-i(\abs{\lambda_+}^2+\abs{\lambda_-}^2)\log\varepsilon} 
    e^{i\abs{\lambda_+}^2\log(\varepsilon+x+t)}
    e^{i\abs{\lambda_-}^2\log(\varepsilon+t-x)},
    \\
    v_\varepsilon(x,t) &= \frac{\lambda_+}{(\varepsilon+x+t)^{1/2}}
    e^{-i(\abs{\kappa_+}^2+\abs{\kappa_-}^2)\log\varepsilon} 
    e^{i\abs{\kappa_+}^2\log(\varepsilon+x+t)}
    e^{i\abs{\kappa_-}^2\log(\varepsilon+t-x)};
  \end{align*}
  \item For $x \le -t$,
  \begin{align*}
    u_\varepsilon(x,t) &= \frac{\kappa_-}{(\varepsilon+t-x)^{1/2}} 
    e^{i\abs{\lambda_-}^2[\log(\varepsilon+t-x) - \log(\varepsilon-x-t)]},
    \\
    v_\varepsilon(x,t) &= \frac{\lambda_-}{(\varepsilon-x-t)^{1/2}} 
    e^{i\abs{\kappa_-}^2[\log(\varepsilon+t-x) - \log(\varepsilon-x-t)]}.
  \end{align*}
\end{enumerate}
Letting $\varepsilon \to 0$, then by the stability assumption, the above should converge to the self-similar solution, but for $t > \abs{x}$ there is no convergence unless $\abs{\lambda_+}^2+\abs{\lambda_-}^2 = 0$ and $\abs{\kappa_+}^2+\abs{\kappa_-}^2 = 0$, in which case the solution is just the trivial one, the zero solution.

\section{Further remarks on the massless case}\label{FinalRemarks}

Assuming $m=0$, we make some further remarks on the behaviour  as $\varepsilon \to 0$ of the approximating solutions $\psi^\varepsilon = (u_\varepsilon,v_\varepsilon)^\intercal$, calculated in the previous section. We take $\kappa_+ = \kappa_- = \lambda_+ = \lambda_- = 1$, so that the data are $f_\varepsilon(x)=g_\varepsilon(x) = (\varepsilon+\abs{x})^{-1/2}$.

Although the solution does not have a limit as $\varepsilon \to 0$, one can nevertheless observe that by restricting $\varepsilon$ to suitable sequences $\varepsilon_n \to 0$, the solution does converge in $C([0,T];L^p([a,b]))$, for $1 \le p < 2$ and any compact interval $[a,b]$, to a valid solution in that space. In this way, one obtains a continuum of possible limiting solutions, depending on the choice of the sequence $\varepsilon_n$. These solutions are all rescalings of each other, via the scaling law \eqref{ScalingLaw}. We now briefly summarize how this works.

Fix $\alpha \in \R$ and choose a positive sequence $\varepsilon_n \to 0$ such that
\[
  e^{-2i\log\varepsilon_n} \longrightarrow e^{i\alpha}
\]
as $n \to \infty$. Letting $\varepsilon$ tend to zero along this sequence $\varepsilon_n$, the solution $(u_\varepsilon,v_\varepsilon)$ above formally converges to $(u,v)$ given by, for $t > 0$,
\begin{enumerate}
  \item For $x > t$,
  \begin{align*}
    u(x,t) &= \frac{1}{(x-t)^{1/2}} e^{i[\log(x+t) - \log(x-t)]},
    \\
    v(x,t) &= \frac{1}{(x+t)^{1/2}} e^{i[\log(x+t) - \log(x-t)]};
  \end{align*}
  \item For $t > \abs{x}$,
  \begin{align*}
    u(x,t) &= \frac{e^{i\alpha}}{(t-x)^{1/2}} 
    e^{i[\log(t-x) + \log(x+t)]},
    \\
    v(x,t) &= \frac{e^{i\alpha}}{(x+t)^{1/2}} 
    e^{i[\log(t-x) + \log(x+t)]};
  \end{align*}
  \item For $x < -t$,
  \begin{align*}
    u(x,t) &= \frac{1}{(t-x)^{1/2}} e^{i[\log(t-x) - \log(-x-t)]},
    \\
    v(x,t) &= \frac{1}{(-x-t)^{1/2}} e^{i[\log(t-x) - \log(-x-t)]}.
  \end{align*}
\end{enumerate}
For $t=0$ we define $u(x,0) = v(x,0) = \abs{x}^{-1/2}$.

In fact, a straightforward analysis reveals that the convergence $(u_\varepsilon,v_\varepsilon) \to (u,v)$ along $\varepsilon=\varepsilon_n$ as $n \to \infty$ is valid in $C([0,T];L^p([a,b]))$ for $1 \le p < 2$ and any $T > 0$.

Moreover, $(u,v)$ makes sense as a solution of the massless Cauchy problem in the sense of distributions on $\R \times (0,T)$, if the nonlinear terms are interpreted in a certain principle value sense, as we show below. Thus, to every point $e^{i\alpha}$ on the unit circle in the complex plane there corresponds a valid solution $(u,v)$ which is obtained as a limit from the family of solutions $(u_\varepsilon,v_\varepsilon)$, and these different solutions all have the same initial data $u(x,0) = v(x,0) = \abs{x}^{-1/2}$. Note that these different solutions are simply rescalings of each other, by the scaling law \eqref{ScalingLaw}.

To end, let us briefly show how the nonlinear terms make sense across the characteristics $x=\pm t$. In wave coordinates
\[
  y = x+t, \qquad s = t-x,
\]
the massless equations become
\begin{align*}
 \partial_y u &= i \abs{v}^2 u,
 \\
 \partial_s v &= i \abs{u}^2 v,
\end{align*}
and we want to check that these are satisfied in a reasonable sense by our functions $u$ and $v$, which are given by (we exclude the quadrant $y,s < 0$, since we consider $t > 0$ only):
\begin{enumerate}
  \item In the region $\Omega_3$ given by $y > 0$ and $s < 0$,
  \begin{align*}
    u &= e^{i\log y} \frac{e^{-i\log(-s)}}{\sqrt{-s}},
    \\
    v &= \frac{e^{i\log y}}{\sqrt{y}} e^{-i\log(-s)};
  \end{align*}
  \item In the region $\Omega_2$ given by $y > 0$ and $s > 0$,
  \begin{align*}
    u(x,t) &= e^{i\alpha} e^{i\log y} \frac{e^{i\log s}}{\sqrt{s}},
    \\
    v(x,t) &= e^{i\alpha} \frac{e^{i\log y}}{\sqrt{y}} 
    e^{i\log s};
  \end{align*}
  \item In the region $\Omega_1$ given by $y < 0$ and $s > 0$,
  \begin{align*}
    u(x,t) &= e^{-i\log(-y)} \frac{e^{i\log s}}{\sqrt{s}},
    \\
    v(x,t) &= \frac{e^{-i\log(-y)}}{\sqrt{-y}} e^{i\log s}.
  \end{align*}
\end{enumerate}

The equations are clearly satisfied in each of the three regions, so what we have to check is what happens across the boundary between $\Omega_1$ and $\Omega_2$, and between $\Omega_2$ and $\Omega_3$. By symmetry, we may limit our attention to the equation
\[
  \partial_y u = i \abs{v}^2 u.
\]

First consider what happens across the boundary between $\Omega_1$ and $\Omega_2$. Then we reduce to the one-dimensional problem
\begin{equation}\label{1dCP}
  f'(y) = \frac{i}{\abs{y}} f(y), \quad y \in \R,
\end{equation}
and we want to check whether the following function satisfies this in a reasonable sense:
\[
  f(y) = \begin{cases}
  e^{-i\log(-y)} &\text{for $y < 0$},
  \\
  e^{i\alpha} e^{i\log y} &\text{for $y > 0$}.
  \end{cases}
\]
This function is in $L^1_{\mathrm{loc}}(\R)$, so it is a distribution and therefore has a derivative on $\R$. It is clear that the equation is satisfied on $(-\infty,0)$ and $(0,\infty)$, but what happens at $y=0$? Fix a test function $\theta \in C_c^\infty(\R)$. Then
\begin{align*}
  \innerprod{f'}{\theta} &= - \innerprod{f}{\theta'}
  \\
  &= - \int_{-\infty}^0 f(y) \theta'(y) \, dy - \int_0^\infty f(y) \theta'(y) \, dy
  \\
  &= - \lim_{\delta \to 0} \left( \int_{-\infty}^{-\delta} f(y) \theta'(y) \, dy
  + \int_\delta^\infty f(y) \theta'(y) \, dy \right).
\end{align*}
Integration by parts then gives
\[
  \innerprod{f'}{\theta} = \lim_{\delta \to 0} \left( \int_{\abs{y} > \delta} \frac{i}{\abs{y}} f(y) \theta(y) \, dy
  + R(\delta) \right),
\]
where
\begin{align*}
  R(\delta) &= e^{i\alpha} e^{i\log \delta}\theta(\delta) - e^{-i\log \delta}\theta(-\delta)
  \\
  &=
  e^{i\alpha} e^{i\log \delta} \left[ \theta(\delta) - \theta(-\delta) \right]
  + \left( e^{i\alpha} e^{i\log \delta} - e^{-i\log \delta} \right) \theta(-\delta).
\end{align*}
Now if we choose $\delta = \delta_n \to 0$ such that
\[
  \lim_{n \to \infty} e^{-2i\log \delta_n} = e^{i\alpha},
\]
then it follows that $\lim_{n \to \infty} R(\delta_n) = 0$, and therefore
\[
  \innerprod{f'}{\theta} = \lim_{n \to \infty} \int_{\abs{y} > \delta_n} \frac{i}{\abs{y}} f(y) \theta(y) \, dy,
\]
which means that the equation in \eqref{1dCP} is satisfied in the sense of distributions on $\R$ when the right hand side is interpreted in a principle value sense for the special intervals $(-\delta_n,\delta_n)$.

Finally, we consider what happens across the boundary between $\Omega_2$ and $\Omega_3$. Then the equation reads
\[
  \partial_y u = \frac{i}{y} u, \quad y > 0, \; s \in \R,
\]
and $u$ is given by, for $y > 0$ and $s \in \R$,
\[
  u = e^{i\log y} F(s)
\]
where
\[
  F(s) = \begin{cases}
  \frac{e^{-i\log(-s)}}{\sqrt{-s}} &\text{for $s < 0$},
  \\
  e^{i\alpha} \frac{e^{i\log s}}{\sqrt{s}} &\text{for $s > 0$}
  \end{cases}
\]
is a locally integrable function on $\R$. It is therefore easy to check that the above equation is satisfied in the sense of distributions in the region $y > 0$ and $s \in \R$.

\bibliographystyle{amsplain}
\bibliography{database}

\end{document}